\documentstyle[12pt]{article}
\def\la{q -q^{-1}}

\def\Si{\Sigma}

\def\de{\delta}

\def\a{a^+}
 
\def\fr{\frac}
\def\q{q^{-1}}

\begin{document}


\vspace{.2cm}
\rightline{GNU-TG-95-08}

\begin{center}
\Large {\bf   On Witten-type deformation of osp(1/2) Algebra}
\\[1cm]
\large W.Chung \\[.3cm]
\normalsize  
Theory Group, Department of Physics, College of Natural Sciences,  \\
\normalsize  Gyeongsang National University,   \\
\normalsize   Jinju, 660-701, Korea
\end{center}

\begin{abstract}
In  this  paper Witten  type    deformation of  osp(1/2)  algbera  is
introduced and  its 
realization and matrix representation are obtained. 
The matrix representation is shown to be possible only when the dimension is
odd.
\end{abstract}

\section{Introduction}

Quantum groups or q-deformation of Lie algebra implies some specific deformation
of classical Lie algebra. From a mathematical point of view, it is a
non-commutative
associative quasi-triangular Hopf algebra. 
The first quantum deformations of the  classical Lie algebra  $su(N)$
treated the elements 
in the Cartan subalgebra on a different
footing [1-4]. For $su(2)$, the deformed algebra is defined as
\begin{equation}
[J_0,J_{\pm}]=\pm J_{\pm},~~[J_+,J_-]=[2H],
\end{equation}
where the q-number is defined as
\begin{equation}
[x]=\fr{q^x-q^{-x}}{\la}.
\end{equation}
This algebra  is sometimes called the transcendental deformation because the
q-number is
given in terms of the transcendental functions.
It is  known   by Woronowicz  [5] and   Witten [6]  that another  type  of
deformation of 
ordinary 
$su(2)$ algebra  is possible .
They treated the generators on a similar footing. And their work enabled
Curtright
and Zachos [7] and Fairlie [8] to find the explicit invertible functionals
that map
su(2) algebra  generators to those of q-deformed algebra . These functionals
are shown to
deform su(2) continuously and reversibly ( except for the special value of
deformation
parameter ) into each of the quantum universal enveloping algebraic structures.
They also applied their method to an obvious two parameter extension which
covers
both Woronowicz's and Witten's form, which is given by
\begin{displaymath}
r HJ_+-r^{-1}J_+H=J_+,
\end{displaymath}
\begin{displaymath}
rJ_-H-r^{-1}HJ_-=J_-,
\end{displaymath}
\begin{equation}
s^{-1}J_+J_--sJ_-J_+=H.
\end{equation}
In this algebra  $s=r^2$ reproduces the Witten's algebra  while $
s=\sqrt{r}$ reproduces
the Woronowicz's algebra .

In this paper we construct the Witten type deformation for the $osp(1/2)$
algebra and
obtain the difference  operator realization  and matrix representation for
this algebra  . 
The triangular deformation 
of
superalgebra is well-known in many papers [9-12].

\section{Witten Type Deformation of OSP(1/2) Algebra}
In this section we discuss the Witten type deformation of $osp(1/2)$ algebra.
Let           us            start                         with           the    
following          
form       
of      
the         
q-deformed  
catesian 
$osp(1/2)$ 
algebra;
\begin{displaymath}
\{V_-, V_+ \}_q=H,
\end{displaymath}
\begin{displaymath}
\{V_-, V_- \}_q=J_-,
\end{displaymath}
\begin{displaymath}
\{V_+, V_+ \}_q=J_+,
\end{displaymath}
\begin{displaymath}
[H,V_+]_q=V_+,
\end{displaymath}
\begin{equation}
[V_-, H]_q =V_-,
\end{equation}
where the qummutator and antiqummutator are defined as
\begin{equation}
[A,B]_q=AB-q BA,~~~\{A,B \}_q=AB+q BA.
\end{equation}
From              the                eq.(4),    especially               the  
second                 
and      
third  
equations,   
one    
can 
obtain 
the
remaining  
commutation                                                              
relations           
of       
the     
Witten type deformation of $osp(1/2)$ 
algebra which is given by
\begin{displaymath}
[J_{\pm}, V_{\pm}]=0,
\end{displaymath}
\begin{displaymath}
[V_-,J_+]_{q^2}=[2]V_+,
\end{displaymath} 
\begin{displaymath}
[J_-,V_+]_{q^2}=[2]V_-,
\end{displaymath} 
\begin{displaymath}
[H,J_+]_{q^2}=[2]J_+,
\end{displaymath} 
\begin{displaymath}
[J_-,H]_{q^2}=[2]J_-,
\end{displaymath} 
\begin{equation}
[J_-,J_+]_{q^4} =[2]^2 (q H +(1-q) V_-V_+),
\end{equation}

where the q-number $[x]$ is defined as $ [x]=\fr{q^x-1}{q-1}$.

If we replace the three bosonic generators by
\begin{displaymath}
H \rightarrow \q [2]H,~~ J_{\pm} \rightarrow \pm q^{-1} [2]^{3/2} J_{\pm}
\end{displaymath}
and decouple the fermionic generators
in order to compare this algbera with its bosonic analogue (3),
 we have the following algebra

\begin{displaymath}
\q H J_+ -qJ_+ H=J_+,
\end{displaymath} 
\begin{displaymath}
\q  J_- H -q HJ_-=J_-,
\end{displaymath} 
\begin{equation}
q^2J_+J_- -q^{-2}J_-J_+= H,
\end{equation}
This correspond to the case that $r=s^2$ in eq.(3), which implies that this
algebra is 
supersymmetric extension of Witten's algebra.

A      glance                  at      the   last
equation   of eq.(6)  
indicates  
that         
the    
qummutator  
of  
the  
two 
even 
generators    (          $         J_+, J_-$)      can        not         be   
written   
in     
terms   
of   
the  
function 
in 
$H$ only                .  However   ,    when          the      deformation    
parameter        
$q$  
goes    
over  
to 
1, 
the 
left              hand              side              of                 the

last     
equation        
of       
eq.(2)     
becomes       
an  
ordinary 
commutator , which is given in terms of the function in $H$ only.

From algebra(4) we can set
\begin{equation}
V_+V_-=F(H),
\end{equation}
where $F(H)$ is a function in $H$ which will be fixed later.
Let $f$ be an inverse function of $F$; then we get
\begin{equation}
H=f(V_+V_-)=\Si_n c_n (V_+V_-)^n.
\end{equation}
From the 4th and 5th relation of (4) we have
\begin{displaymath}
[H,V_+]_q
\end{displaymath}
\begin{displaymath}
=[\Si_n c_n (V_+V_-)^n, V_+]_q
\end{displaymath}
\begin{displaymath}
=\Si_n c_n [ V_+(V_+V_-)^n-q V_+(V_+V_-)^n]
\end{displaymath}
\begin{equation}
=V_+[f(V_-V_+)-qf(V_+V_-)]=V_+.
\end{equation}
Thus we have
\begin{displaymath}
f(V_-V_+)=qf(V_+V_-)+1,
\end{displaymath}
which implies
\begin{equation}
V_-V_+=F(qH+1).
\end{equation}
Inserting (8) and (10) into the 1st relation of (4) gives
\begin{equation}
F(1+qH)+qF(H)=H.
\end{equation}
We can easily obtain the Casimir operator for the algebra (4) as  follows;
\begin{equation}
C=V_+V_--\fr{1}{2q}(H-\fr{1}{1+q}).
\end{equation}

Let us  denote by $|c,m>$ a  simultaneous eigenvector of the  commuting
hermitian operator 
$H$
and $C$, where we have
\begin{equation}
H|c,m>=m|c,m>,~~C|c,m>=c|c,m>.
\end{equation}
From the algebra (4) and the fact that $V_-$ is a hermitian conjugate of
$V_+$, we see
 that
\begin{equation}
V_+|c,m>=a(c,qm+1)|c,qm+1>,
\end{equation}
\begin{equation}
V_-|c,m>=a(c,m)|c,q^{-1}m -q^{-1}>.
\end{equation}
Using the Casimir operator and algebra (4) we have
\begin{equation}
V_-|c,m>=\sqrt{\fr{j+m}{2q}}|c,q^{-1}(m-1)>,
\end{equation}
\begin{equation}
V_+|c,m>=\sqrt{\fr{j+qm+1}{2q}}|c,qm+1)>,
\end{equation}
where we set $ m\geq -j$.
If we define
\begin{displaymath}
|j,n>=|c,-q^n j +[n]>,~~[n]=\fr{1-q^n}{1-q}
\end{displaymath}
then we have a simpler representation of algebra (4) as follows;
\begin{displaymath}
H|j,n>=([n]-q^n j)|j,n>,~~(n=0,1,2,\cdots)
\end{displaymath}
\begin{displaymath}
V_-|j,n>=\sqrt{\fr{1+(1-q)j}{2q}}\sqrt{[n]}|j,n-1>,
\end{displaymath}
\begin{equation}
V_+|j,n>=\sqrt{\fr{1+(1-q)j}{2q}}\sqrt{[n+1]}|j,n+1>.
\end{equation}
So the representation is infinite dimensional and is bouned from below.

\section{Realization}
In  this section  we discuss  the  realization of  Witten type  deformation
of $osp(1/2)$ 
algebra.
In doing so we need the new type of the deformed oscillator algebra.
Now let us consider the following deformed oscillator algebra
\begin{equation}
[a,\a]_q=1,~~[N,\a]_q=\a,~~[a,N]_q=a
\end{equation}
where $a$ is assumed to be a hermitian conjugate of $\a$ and the number
operator $N$ is
assumed to be hermitian. Then the third relation of eq.(20) is not necessary
any more
because it is obtained by taking the complex conjugate of the second
relation of eq.(20).
The first relation itself takes self-hermitian form. If we replace the
qummutators with
the ordinary commutators in the second and the third relation of eq.(20),
this algebra
reduces to the standard q-deformed boson algebra.
From the definition (20), we have
\begin{equation}
\a a =N,~~ a \a =1+qN .
\end{equation}
If we assume that there exists a unique vacuum state $|0>$ satisfying
\begin{equation}
a|0>=0,
\end{equation}
we have, from the relation (21),
\begin{equation}
N|0>=0.
\end{equation}
Then the higher states can be obtained by acting the creation operator $\a$
 on the vacuum state successively. Let the state $|n>$ be proportional to
$(\a)^n|0>$.
 Then we have
\begin{equation}
N|n>=[n]|n>,~~n=0,1,2,\cdots.
\end{equation}
So the eigenvalues of the number operator $N$ become a sequence of the
q-numbers.
Using the second and third relation of eq.(20) , we get
\begin{equation}
a|n>=\sqrt{[n]}|n-1>,~~\a |n>=\sqrt{[n+1]}|n+1>.
\end{equation}

From this  new deformed  boson algebra we  can obtain  the realization for
the $osp(1/2)$ 
algebra as follows;
\begin{displaymath}
v_-=\fr{1}{\sqrt{1+q}}a
\end{displaymath}
\begin{displaymath}
v_+=\fr{1}{\sqrt{1+q}}\a
\end{displaymath}
\begin{equation}
H=\fr{1}{1+q}( a\a + q \a a).
\end{equation}
Now we can obtain the difference operator realization for the $osp(1/2)$
algebra.
The coordinate realization of algebra (20) is given by
\begin{equation}
\a=x,~~a=D,~~N=xD,
\end{equation}
where the q-derivative $D$ is defined as
\begin{equation}
Df(x)=\fr{f(qx)-f(x)}{x(q-1)}
\end{equation}
Using this  relation,  we  can write the difference  operator  realization
for the Witten 
type deformation of $osp(1/2)$ 
algebra as follows;
\begin{equation}
v_-=\fr{1}{\sqrt{1+q}}D,~~v_+=\fr{1}{\sqrt{1+q}} x,~~H=\fr{1}{1+q}(1+2q x D)
\end{equation}

\section{Matrix Representation}
Now     we     will      obtain           the   matrix        representation 
of       
the    
algebra     
(4)  
by 
assuming 
that  $H$ is diagonal and $V_+(V_-)$ are  super (sub) diagonal;
\begin{displaymath}
(H)_{ij}=h_i\de_{ij},
\end{displaymath}
\begin{displaymath}
(V_+)_{ij}=v_i \de_{i+1,j},
\end{displaymath}
\begin{equation}
(V_-)_{ij}=v_j \de_{i,j+1},
\end{equation}
where     the         last        relation     is        obtained       from

the        
second 
relation  
by    
using 
the 
fact          that        $V_-$         is           hermitian     conjugate

to           
$V_+$.   
Then        
the  
forth       
relation 
of eq.(4) determines $h_i$ through the recurrence relation;
\begin{equation}
h_i-qh_{i+1}=1,~~(i=1,2,\cdots, n-1).
\end{equation}
The relation (31) is easily solved and the solution is given by
\begin{equation}
h_i=q^{n-i} h_n +\fr {q^{-1}(q^{n-i}-1)}{1-q^{-1}},
\end{equation}
where $h_n$ will be determined by using the another relations of eq.(4).
Using the first relation of eq.(4), we get
\begin{equation}
v_{i-1}^2 +q v_i^2 =h_i,~~(i=1,2,\cdots, n);~~v_0=0~and ~v_n=0.
\end{equation}
Using           the    eq.(32)       and            checking              the   
consistency    
of       
eq.(33),  
we    
obtain  
the 
following identity;
\begin{equation}
\Si_{k=1}^n (-q^{-1})^{n-k}h_k =0.
\end{equation}
This       relation          and    eq.(32)       fix    the         concrete 
form     
of  
the     
value  
of  
$h_n$. 
Inseting eq.(32) into the eq.(34) produces
\begin{equation}
h_n \fr {(-)^n-1}{2}
=\fr{1}{1-q} \fr {(-)^n-1}{2}
+\fr{1}{1-q}\fr{q+(-q^{-1})^{n-1}}{1+q}.
\end{equation}
For the          case that     $n$     is     even,    the      eq.(35)
does    
not    
give  
the  
solution  
for 
$h_n$; so we assume that  $n$ is odd. Then we have
\begin{equation}
h_n =\fr{1}{1-q} +\fr{q+q^{-n+1}}{(q-1)(q+1)}.
\end{equation}
Inserting eq.(36) into eq.(32) we obtain
\begin{equation}
h_i =\fr{1}{1-q}\{1-\fr{q}{1+q}(q^n+1)q^{-i}\}.
\end{equation}
Using eq.(33)                we   can          perform   the
partial    
summations   
and     
obtain   
a  
closed 
expression for $v_i$ given by;
\begin{equation}
v_i =\sqrt { \fr{1}{1-q^2}[1-(-q^{-1})^i+\fr{(-)^i-1}{2}(q^n+1)q^{-i}]}.
\end{equation}

For $n=3$ case we can obtain the explicit matrix representation as follows;
\begin{displaymath}
H=\left(\matrix{ q& 0&0 \cr 0&1-q^{-1} &0\cr 0&0& -q^{-2}\cr}\right)
\end{displaymath}
\begin{displaymath}
V_+=\left(\matrix{ 0& 1&0 \cr 0&0& iq^{-1}\cr 0&0& 0\cr}\right)
\end{displaymath}
\begin{equation}
V_-=\left(\matrix{ 0& 0&0 \cr 1&0 &0\cr 0&-iq^{-1}&0\cr}\right),
\end{equation}
where we assumed that $q$ is real.

\section*{Acknowledgement}
This                   paper                was
supported         by  
                 
NON             
DIRECTED     
RESEARCH     
FUND,  
Korea   
Research 
Foundation,   1994   and   by   the   KOSEF    through   C.T.P.   at   Seoul  
National 
University.   And   the   present   studies    were   supported   by   Basic  
Science 
Research Program, Ministry of Education, 1995 (BSRI-95-2413).

\vfill\eject

\end{document}